\documentclass[11pt]{article} 
\usepackage{amsmath}
\usepackage{amsfonts}
\usepackage{graphicx}
\usepackage{color}
\usepackage[latin1]{inputenc}
\usepackage[T1]{fontenc}
\usepackage[french]{babel}
\makeatletter
\def\@sect#1#2#3#4#5#6[#7]#8{%
  \ifnum #2>\c@secnumdepth
    \let\@svsec\@empty
  \else
    \refstepcounter{#1}%
    \protected@edef\@svsec{\@seccntformat{#1}\relax}%
  \fi
  \@tempskipa #5\relax
  \ifdim \@tempskipa>\z@
    \begingroup
      #6{%
        \@hangfrom{\hskip #3\relax\@svsec}%
          \interlinepenalty \@M #8\@@par}%
    \endgroup
    \csname #1mark\endcsname{#7}%
    \addcontentsline{toc}{#1}{%
      \ifnum #2>\c@secnumdepth \else
        \protect\numberline{\csname the#1\endcsname.}%
      \fi
      #7}%
  \else
    \def\@svsechd{%
      #6{\hskip #3\relax
      \@svsec #8}%
      \csname #1mark\endcsname{#7}%
      \addcontentsline{toc}{#1}{%
        \ifnum #2>\c@secnumdepth \else
          \protect\numberline{\csname the#1\endcsname.}%
        \fi
        #7}}%
  \fi
  \@xsect{#5}}
\def\@seccntformat#1{\csname the#1\endcsname.\quad}
\makeatother
\newtheorem{theo}[equation]{Th\'eor\`eme}
\newtheorem{cor}[equation]{Corollaire}
\newtheorem{question}[equation]{Question}

\makeatletter
\renewcommand\theequation{\thesection.\arabic{equation}}
\@addtoreset{equation}{section}
\makeatother
\newcommand{\carrenoir}{\rule{0.5em}{0.5em}}
\makeatletter
\newenvironment{demo}[1][\@empty]{\textbf{D\'emonstration~%
\ifx\@empty#1:\else #1~:\fi~}}
{\hfill\carrenoir\nolinebreak\vspace{2mm}}
\makeatother
\newcommand{\oper}[2]{\newcommand{#1}{\mathop{\mathrm{#2}}\nolimits} }
\oper{\Vol}{Vol}
\newcommand{\de}{\mathrm{ d }}
\oper{\dimension}{dim}
{\endtrivlist}%
\newenvironment{remarque}{
\refstepcounter{equation}\trivlist%
\item[\hskip \labelsep{\bfseries Remarque \theequation.\ }]}%

\DeclareSymbolFont{greek}{OML}{ptmcm}{m}{it}
\DeclareMathSymbol{\codiff}{\mathord}{greek}{"0E}

\title{Minoration conforme du spectre du laplacien de Hodge-de~Rham}
\author{Pierre Jammes}
\date{}
\begin{document}
\maketitle
{\small 
\textsc{Résumé.---}
Soit $M^n$ une variété compacte de dimension $n\geq3$. Pour toute
classe conforme $C$ de métriques riemanniennes sur $M$, on pose 
$\mu_k^c(M,C)=\inf_{g\in C}\mu_{\left[\frac n2\right],k}(M,g)
\Vol(M,g)^{\frac2n}$, où $\mu_{p,k}(M,g)$ est la $k$-ième valeur propre 
du laplacien de Hodge-de~Rham agissant sur les $p$-formes coexactes. 
On démontre que $0<\mu_k^c(M,C)\leq\mu_k^c(S^n,[g_\textrm{can}])\leq 
k^{\frac2n}\mu_1^c(S^n,[g_\textrm{can}])$. On montre aussi que
si $n=0,2,3\textrm{ mod }4$ et si $g$ est une metrique lisse telle que 
$\mu_{\left[\frac n2\right],1}(M,g)\Vol(M,g)^{\frac2n}=\mu_1^c(M,[g])$,
alors il existe une forme propre non nulle de degré $\left[\frac{n-1}2\right]$
et de valeur propre $\mu_1^c(M,[g])$ qui est de longueur constante. 
En conséquence, il n'existe pas de métrique extrémale lisse si $n=4$ et
que la caractéristique d'Euler de $M$ est non nulle.

Mots-clefs : formes différentielles, laplacien, métriques conformes,
inégalités de Sobolev.

\medskip
\textsc{Abstract.---}
Let $M^n$ be a $n$-dimensional compact manifold, with $n\geq3$. For any 
conformal class $C$ of riemannian metrics on $M$, we
set $\mu_k^c(M,C)=\inf_{g\in C}\mu_{\left[\frac n2\right],k}(M,g)
\Vol(M,g)^{\frac2n}$, where $\mu_{p,k}(M,g)$ is the $k$-th 
eigenvalue of the Hodge laplacian acting on coexact $p$-forms.  
We prove that $0<\mu_k^c(M,C)\leq \mu_k^c(S^n,[g_\textrm{can}])\leq
 k^{\frac2n}\mu_1^c(S^n,[g_\textrm{can}])$. We also prove that if $g$
is a smooth metric such that $\mu_{\left[\frac n2\right],1}(M,g)
\Vol(M,g)^{\frac2n}=\mu_1^c(M,[g])$, and $n=0,2,3\textrm{ mod }4$, then 
there is a non-zero corresponding eigenform of degree 
$\left[\frac{n-1}2\right]$ with constant length. As a
corollary, on a 4-manifold with non vanishing Euler characteristic,
there is no such smooth extremal metric. 


}

\section{Introduction}
Soit $(M,g)$ une variété riemannienne compacte connexe orientable
de dimension $n$. 
On considère le laplacien de Hodge-de~Rham agissant sur l'espace
$\Omega^p(M)$ des $p$-formes différentielles sur $M$ et défini par 
$\Delta=\de\codiff+\codiff\de$ où la codifférentielle $\codiff$ est
l'adjoint $L^2$ de $\de$. Le spectre de cet 
op\'erateur forme un ensemble discret de nombres positifs ou nuls qu'on
notera
$$0=\lambda_{p,0}(M,g)<\lambda_{p,1}(M,g)\leq\lambda_{p,2}(M,g)\leq\dots,$$
où la multiplicité de $\lambda_{p,0}(M,g)$ est le $p$-ième nombre de
Betti de $M$, les autres valeurs propres étant répétée s'il y a 
multiplicité.

L'étude du laplacien agissant sur les fonctions, c'est-à-dire du cas $p=0$, 
montre qu'à volume fixé, on peut choisir la métrique de manière à rendre 
les valeurs propres arbitrairement grande (voir \cite{cd94}). mais
que ce n'est pas le cas si on impose à la métrique d'appartenir à une
classe conforme donnée, c'est-à-dire que si $C$ est une classe
conforme de métrique riemannienne sur $M$, alors
\begin{equation}\label{intro:eq1}
\sup_{g\in C}\lambda_{0,k}(M,g)\Vol(M,g)^{\frac2n}<+\infty,
\textrm{ pour tout }k.
\end{equation}
Cette majoration a été démontrée par A.~El~Soufi et S.~Ilias pour $k=1$
(\cite{esi86}), et généralisée par N.~Korevaar dans \cite{ko93}. On
peut par ailleurs facilement montrer que $\inf_{g\in C}\lambda_{0,k}(M,g)
\Vol(M,g)^{\frac2n}=0$.

Le laplacien agissant sur les formes différentielles de degré quelconque
a été moins étudié. B.~Colbois et A.~El~Soufi ont cependant montré
récemment que l'inégalité (\ref{intro:eq1}) ne se généralise pas:
\begin{theo}[\cite{ces06}]\label{intro:th1}
Si $M$ est une variété riemannienne compacte et $C$ une classe conforme sur 
$M$, alors
\begin{equation}
\sup_{g\in C}\inf_{2\leq p\leq n-2}\lambda_{p,1}(M,g)\Vol(M,g)^{\frac2n}=
+\infty.
\end{equation}
\end{theo}
Comme la différentielle $\de$ commute avec le laplacien, le
spectre des $1$-formes contient le spectre des fonctions, on peut
donc majorer $\lambda_{1,k}$ comme dans (\ref{intro:eq1}), et il
en va de même pour les $n-1$ formes par dualité de Hodge.

Le théorème \ref{intro:th1} conduit naturellement à se demander si on
peut faire tendre les valeurs propres du laplacien de Hodge-de~Rham
vers zéro dans une classe conforme et à volume fixé, comme c'est le cas
pour les fonctions:
\begin{question}[\cite{co04}]\label{intro:q1}
Étant donnés $2\leq p\leq n-2$ et $k\geq1$, a-t-on
\begin{equation}
\inf_{g\in C}\lambda_{p,k}(M,g)\Vol(M,g)^{\frac2n}=0\ ?
\end{equation}
\end{question}
 On montre ---~entre autres choses~--- dans \cite{ja06} que la réponse est 
positive pour tout $p$ et $k$
si la dimension $n$ est impaire, et pour tout $k$ et tout $p\neq\frac n2$
si $n$ est pair. Plus précisément, si on se restreint aux formes 
différentielles coexactes, la technique utilisée échoue en degré $\frac n2$ et
$\frac n2-1$ si la dimension est paire, et en degré $\left[\frac n2\right]$
si la dimension est impaire, la question \ref{intro:q1} restant ouverte
dans ces cas.

Nous allons ici achever de répondre à la question \ref{intro:q1} en montrant
que dans les cas restants, on peut en fait minorer uniformément le spectre
des formes coexactes, le volume et la classe conforme étant fixés, par
une constante strictement positive. La réponse à la question \ref{intro:q1}
dépend donc en général du degré $p$.

Nous nous appuieront sur l'existence d'inégalités de Sobolev pour les
formes différentielles, et en particulier de l'inégalité suivante,
dont on peut par exemple trouver la démonstration dans \cite{gt05}:
\begin{theo}
Soit $(M,g)$ une variété compacte de dimension $n$, et $r,s\in]1,+\infty[$
deux réels tels que $1/s-1/r=1/n$.
Il existe une constante $c(M,g,r,s)>0$ telle que pour toute $p$-forme 
$\theta\in\Omega^p(M,g)$, on a
\begin{equation}\label{intro:eq2}
\inf_{\de\zeta=0}\|\theta-\zeta\|_{L^r}\leq c\|\de\theta\|_{L^s}.
\end{equation} 
\end{theo}
\begin{remarque}
La démonstration donnée dans \cite{gt05} étant en grande partie
non constructive, elle ne fournit aucune estimée de la constante $c$.
\end{remarque}
Dans le cas où $r=\frac np$ et $s=\frac n{p+1}$, les normes 
$\|\theta-\zeta\|_r$ et $\|\de\theta\|_s$ sont conformément invariantes,
et donc la constante $c$ optimale dans l'inégalité (\ref{intro:eq2})
en restriction aux $p$-formes est un invariant conforme, qu'on notera
\begin{equation}
K_p(M,[g])=\sup_{\theta\in\Omega^p(M)}\inf_{\de\zeta=0}
\frac{\|\theta-\zeta\|_{\frac np}}{\|\de\theta\|_{\frac n{p+1}}},
\end{equation}
où $[g]=\{h^2g, h\in C^\infty(M), h>0\}$ désigne la classe conforme de $g$.

On va donner une minoration uniforme du spectre du laplacien de 
Hodge-de~Rham en fonction de la constante $K_p$. On notera 
$0<\mu_{p,1}(M,g)\leq\mu_{p,2}(M,g)\leq\ldots$ le spectre du laplacien
en restriction aux $p$-formes coexactes. Rappelons que la spectre
non nul $(\lambda_{p,k}(M,g))_{k\geq1}$ du laplacien sur les $p$-formes
est la réunion des spectres $(\mu_{p,k}(M,g))$ et $(\mu_{p-1,k}(M,g))$.
\begin{theo}\label{intro:th2}
Soit $M^n$ une variété compacte de dimension $n\geq3$, $C$ une classe
conforme de métriques sur $M$. Pour
toute métrique $g\in C$, on a
$$\mu_{\left[\frac n2\right],1}(M,g)\Vol(M,g)^{\frac2n}\geq 
K_{\left[\frac n2\right]}(M,C)^{-2}.$$
Si $n$ est pair, cette inégalité est optimale. 
\end{theo}
\begin{remarque} La théorie de Hodge nous dit que si $n$ est pair, alors
$\mu_{\frac n2-1,i}(M,g)=\mu_{\frac n2,i}(M,g)$ 
pour tout $i\geq1$. Comme on l'a déjà remarqué, pour tous les autres 
degrés $p$, on sait faire tendre
$\mu_{p,i}(M,g)$ vers zéro à volume fixé dans la classe conforme $C$. 
On a donc bien complètement répondu à la question \ref{intro:q1}.
\end{remarque}

Dans \cite{ces03} B.~Colbois et A.~El~Soufi définissent un «~spectre
conforme~» pour les fonctions par $\lambda_k^c(M,C)=\sup_{g\in C}
\lambda_{0,k}(M,g)\Vol(M,g)^{\frac2n}$. On peut définir de la même manière 
un spectre conforme pour les formes différentielles par 
\begin{equation}
\mu_k^c(M,C)=\inf_{g\in C}\mu_{\left[\frac n2\right],k}(M,g)
\Vol(M,g)^{\frac2n}.
\end{equation}

 La théorème \ref{intro:th2} soulève naturellement la question de 
l'existence de métriques lisses réalisant $\mu_k^c(M,C)$, en particulier 
pour $k=1$. On peut établir le critère suivant :
\begin{theo}\label{intro:th3}
Soit $g$ une métrique lisse de volume~1 sur $M$ telle que 
$\mu_{\left[\frac n2\right],1}(M,g)=\mu_1^c(M,[g])$. Si $n=3\textrm{ mod }4$,
alors il existe une $\left[\frac n2\right]$-forme propre coexacte non nulle
de valeur propre $\mu_1^c(M,[g])$ et
de longueur constante. Si $n$ est pair, toutes les formes propres coexacte 
de degré $\frac n2-1$ et de valeur propre $\mu_1^c(M,[g])$ sont de 
longueur constante.
\end{theo}
\begin{remarque}
On sait que pour les fonctions, la métrique canonique de la sphère
maximise la première valeur propre dans sa classe conforme (\cite{esi86}).
En revanche les premières formes propres des formes
différentielles pour cette métrique ne sont pas de longueur constante
(voir l'appendice de \cite{GM75}). En général, les métriques extrémales ne
sont donc pas les mêmes pour les formes et les fonctions.
\end{remarque}
On peut déduire du théorème \ref{intro:th3} qu'en dimension~4, 
il y a une obstruction topologique à l'existence de métriques
réalisant l'infimum $\mu_1^c(M,C)$:

\begin{cor}\label{intro:cor}
 Si $M$ est une variété compacte de dimension~4 et de
caractéristique d'Euler non nulle, il n'existe pas de métrique régulière
$g$ de volume~1 telle que $\mu_{1,1}(M,g)=\mu_1^c(M,[g])$.
\end{cor}

Enfin, on va donner une majoration du spectre conforme 
$(\mu_k^c(M,C))_{k\geq1}$:
\begin{theo}\label{intro:th4}
Pour toute classe conforme $C$ sur la variété compacte $M$ et tout $k\geq1$, 
on a
$$\mu_k^c(M,C)\leq\mu_k^c(S^n,C_\textrm{\emph{can}})\leq k^{\frac2n}
\mu_1^c(S^n,C_\textrm{\emph{can}}),$$
où $C_\textrm{\emph{can}}$ est la classe conforme de la métrique canonique 
de la sphère.
\end{theo}

Il existe une minoration semblable à celle du théorème \ref{intro:th2}
pour l'opérateur de Dirac (voir \cite{lo86} et \cite{am03}). L'étude
de métriques extrémales pour cet opérateur a été développée par
B.~Ammann dans \cite{am03b}, où il montre entre autres l'optimalité
de l'équivalent spinoriel du théorème \ref{intro:th2}. Cependant,
les techniques utilisées ne sont en général pas immédiatement
transposables au cas des formes différentielles; en particulier l'évolution 
de l'opérateur de Dirac au cours d'une déformation conforme s'exprime
beaucoup plus simplement que celle du laplacien de Hodge-de~Rham.

Les théorèmes \ref{intro:th2} et \ref{intro:th3} et le corollaire
\ref{intro:cor} seront démontrés dans
la section \ref{min} et le théorème \ref{intro:th4} dans la section
\ref{maj}.

\section{Extremum du spectre}\label{min}
On va voir ici comment l'inégalité de Sobolev (\ref{intro:eq2})
permet de minorer le spectre.

\begin{demo}[du théorème \ref{intro:th2}] 
Pour commencer, on considérera le cas où la dimension $n$ est paire.

Soit $\theta$ une forme propre coexacte de degré $\frac n2$, de valeur 
propre $\lambda$. On sait que $\|\de\theta\|_2^2=\lambda\|\theta\|_2^2$,
et en supposant que la métrique $g$ sur $M$ est de volume~1, on peut en outre 
écrire $\|\de\theta\|_2^2\geq\|\de\theta\|_{\frac{2n}{n+2}}^2$ (inégalité
de Hölder). L'inégalité
de Sobolev (\ref{intro:eq2}), en conjonction avec le fait $\theta$
minimise la norme $L^2$ parmi les primitives de $\de\theta$, assure alors que
\begin{eqnarray}
\|\theta\|_2^2&\leq&K_{\frac n2}(M,[g])^2\|\de\theta\|_{\frac{2n}{n+2}}^2
\leq K_{\frac n2}(M,[g])^2\|\de\theta\|_2^2\nonumber\\
&\leq& K_{\frac n2}(M,[g])^2 \lambda\|\theta\|_2^2.
\end{eqnarray}
On a donc bien $\lambda\geq K_{\frac n2}(M,[g])^{-2}$.

En dimension paire, on peut aussi montrer que cette minoration est optimale.
Soit $\varepsilon>0$, et $\theta$ une $\frac n2$-forme coexacte telle que
$\|\de\theta\|_{\frac{2n}{n+2}}=1$ et
\begin{equation}
K_{\frac n2}(M,[g])-\varepsilon\leq\|\theta\|_2\leq K_{\frac n2}(M,[g]).
\end{equation}

Supposons dans un premier
temps que la forme $\de\theta$ soit partout non nulle. Sa longueur
est donc partout non nulle et on peut définir une nouvelle métrique
par $g_h=h^2g$ avec $h=|\de\theta|^{\frac2{n+2}}$; on notera
$|\cdot|_h$ et $\|\cdot\|_{p,h}$ les normes ponctuelle et $L^p$ pour
cette métrique. Le volume pour la nouvelle métrique est 
$\Vol(M,g_h)=\int_M h^n\de v_g=
\|\de\theta\|_{\frac{2n}{n+2}}^{\frac{2n}{n+2}}=1$, et surtout la longueur
de la forme $\de\theta$ pour la métrique $g_h$ est constante: en effet, 
$|\de\theta|_h=h^{-(\frac n2+1)}|\de\theta|=1$. La forme $\theta$
étant de degré $\frac n2$, elle est coexacte quel que soit le choix de
la métrique dans $[g]$,
on peut donc facilement majorer $\mu_{\frac n2,1}(M,g_h)$ à l'aide de 
son quotient de Rayleigh:
\begin{equation}\label{min:eq1}
\mu_{\frac n2,1}(M,g_h)\leq\frac{\|\de\theta\|_{2,h}^2}{\|\theta\|_{2,h}^2}=
\frac{\|\de\theta\|_{\frac{2n}{n+2},h}^2}{\|\theta\|_{2,h}^2}\leq
(K_{\frac n2}(M,[g])-\varepsilon)^{-2}.
\end{equation}
En passant à la limite quand $\varepsilon\to0$, on obtient bien que
$\mu_1^c(M,[g])=K_{\frac n2}(M,[g])^{-2}$.

Si la forme $\de\theta$ n'est pas partout non nulle, on doit modifier 
légèrement la démonstration. On se donne un réel $\delta>0$ et on pose
$h=(|\de\theta|+\delta)^{\frac2{n+2}}$. On a alors
$|\de\theta|_h=\frac{|\de\theta|}{|\de\theta|+\delta}\leq1$, et donc
$\|\de\theta\|_{2,h}\leq\|\de\theta\|_{\frac{2n}{n+2},h}$. Comme
$\Vol(M,g_h)\to1$ quand $\delta\to0$, on retrouve l'inégalité (\ref{min:eq1})
en passant à la limite.

Considérons maintenant la situation où $n$ est impaire. Soit $\theta$
une forme propre coexacte de degré $\left[\frac n2\right]$ et de
valeur propre $\lambda$. Par définition de $K_p$, pour tout $\varepsilon>0$,
il existe une forme $\theta'$ de degré $\left[\frac n2\right]$ telle que
$\|\theta'\|_{\frac{2n}{n-1}}\leq(K_{\frac n2}(M,[g])+\varepsilon)
\|\de\theta\|_{\frac{2n}{n+1}}$. On a par conséquent
\begin{eqnarray}
\|\theta\|_2^2&\leq&\|\theta'\|_2^2\leq\|\theta'\|_{\frac{2n}{n-1}}^2
\leq(K_{\frac n2}(M,[g])+\varepsilon)^2\|\de\theta\|_{\frac{2n}{n+1}}^2
\nonumber\\
&\leq&(K_{\frac n2}(M,[g])+\varepsilon)^2\|\de\theta\|_2^2.
\end{eqnarray}
En passant à la limite quand $\varepsilon\to0$, on obtient la même conclusion
qu'en dimension paire.
\end{demo}

On va maintenant établir le critère du théorème \ref{intro:th3} pour
l'existence de métriques extrémales.

\begin{demo}[du théorème \ref{intro:th3}]
 On se donne une métrique $g$ sur $M$ de volume~1 telle que 
$\mu_{\left[\frac n2\right],1}(M,g)=\mu_1^c(M,[g])$. On va ici encore donner 
deux démonstrations distinctes selon la parité de la dimension. 

Dans le cas où la dimension $n$ est paire, il est commode de raisonner
sur la longueur de formes propres exactes de degré $\frac n2+1$ et
de conclure par dualité de Hodge: si $\omega$ est une forme propre
coexacte de degré $\frac n2-1$, alors $*\omega$ est une forme propre
exacte de degré $\frac n2+1$ de même longueur et de même valeur propre.

Considérons donc une forme propre exacte $\omega$ de degré $\frac n2+1$ et
de valeur propre $\mu_{\left[\frac n2\right],1}(M,g)$, et notons $\varphi$
la forme propre coexacte de degré $\frac n2$ telle que $\de\varphi=\omega$.
On va étudier les variations du quotient de Rayleigh de $\varphi$ pour
des petites déformations conformes de la métrique. On se donne une
fonction $u\in C^\infty(M)$, et on définit pour tout $\varepsilon>0$ petit 
la métrique $g_\varepsilon=(1+\varepsilon u)^2g$. On peut imposer les 
conditions de normalisation que la volume reste constant au 1\ier~ordre
par rapport à $\varepsilon$ c'est-à-dire que $\int_M u\de v_g=0$,
et que $\|\varphi\|_2=1$.

Comme $\varphi$ est de degré $\frac n2$, elle reste coexacte et de norme~1
quel que soit~$\varepsilon$. Le développement de son quotient de Rayleigh
est
\begin{eqnarray}
R_{g_\varepsilon}(\varphi)&=&\|\omega\|_{2,\varepsilon}^2=
\int_M|\omega|_{g_\varepsilon}^2\de v_{g_\varepsilon}=
\int_M|\omega|^2(1+\varepsilon u)^{-2}\de v_g\nonumber\\
&=&\int_M|\omega|^2\de v_g-2\varepsilon\int_M u|\omega|^2\de v_g
+o(\varepsilon).
\end{eqnarray}
Comme la métrique $g$ minimise $\mu_{\left[\frac n2\right],1}(M,g)$
dans sa classe conforme, le terme d'ordre~1 du développement précédent est
nul, c'est-à-dire que $\int_M u|\omega|^2\de v_g=0$, quelle que soit
la fonction $u$ d'intégrale nulle, ce qui n'est possible que si $|\omega|$
est constant.

Traitons maintenant le cas où $n=3\textrm{ mod }4$. On va en fait
raisonner sur l'opérateur $(*\de)$ agissant sur l'espace des formes
coexactes de degré $\left[\frac n2\right]$. C'est un opérateur autoadjoint
qui est une racine carrée du laplacien et qui laisse stable les espaces 
propres de $\Delta$. Pour chaque valeur propre $\mu$ du laplacien sur cet
espace, l'espace propre associé est donc un espace propre de $(*\de)$ ou 
se décompose en la somme orthogonale de deux espaces propres de 
$(*\de)$ de valeurs propres $\sqrt\mu$ et $-\sqrt\mu$. Une forme
propre $\varphi$ de $(*\de)$ a la particularité de vérifier 
$|\de\varphi|=\sqrt\mu|\varphi|$ en tout point. On note $\omega$
une forme propre coexacte du laplacien de norme~1,
de valeur propre $\mu_{\left[\frac n2\right],1}(M,g)$ et qui est aussi
forme propre de $(*d)$,
et on va calculer le développement au 1\ier~ordre du quotient de Rayleigh 
(pour le laplacien) de $\omega$ par rapport à $\varepsilon$. 

Remarquons d'abord qu'au 1\ier~ordre, la norme 
$\|\codiff\omega\|_{g_\varepsilon}^2$ est nulle: en effet, 
$\|\codiff\omega\|_{g_\varepsilon}=\|\de*_{g_\varepsilon}\omega
\|_{g_\varepsilon}=\|\de[(1+\varepsilon u)*_g\omega]\|_{g_\varepsilon}=
\varepsilon\|(1+\varepsilon u)^n\de u\wedge*\omega\|_g$. On peut alors écrire 
\begin{eqnarray}
R_{g_\varepsilon}(\omega)&=&\frac{\|\de\omega\|_{g_\varepsilon}^2}
{\|\omega\|_{g_\varepsilon}^2}+o(\varepsilon)\nonumber\\
&=&\frac{\int_M|\de\omega|_g^2(1+\varepsilon u)^{-1}\de v_g}
{\int_M|\omega|_g^2(1+\varepsilon u)\de v_g}+o(\varepsilon)\nonumber\\
&=&\|\de\omega\|_g^2-\varepsilon\int_Mu(|\de\omega|^2+
\|\de\omega\|^2|\omega|^2)\de v_g.
\end{eqnarray}
 On conclut comme dans le cas de la dimension paire en utilisant le fait
que $|\de\omega|^2=\|\de\omega\|^2|\omega|^2$ en tout point.
\end{demo}
\begin{remarque}
Si $n=1\textrm{ mod }4$, la démonstration précédente n'est pas valide car
l'opérateur $(*\de)$ n'est pas autoadjoint.
\end{remarque}

\begin{demo}[du corollaire \ref{intro:cor}]
Il suffit de remarquer qu'en dimension~4, L'existence d'une métrique
extrémale implique l'existence d'une $1$-forme partout non nulle, et
donc que la caractéristique d'Euler est nulle.
\end{demo}

\section{Majoration du spectre conforme}\label{maj}

Pour finir, nous allons démontrer les majorations du spectre conforme
énoncées dans l'introduction.

\begin{demo}[du théorème \ref{intro:th4}]%
Pour démontrer l'inégalité $\mu_k^c(M,C)\leq\mu_k^c(S^n,C_\textrm{can})$,
on va montrer qu'on peu faire tendre le spectre et le volume de $M$ vers 
celui d'une sphère dans la classe conforme $C$. On s'appuiera sur
un théorème de convergence de spectre obtenu par C.~Anné et B.~Colbois 
dans \cite{ac95} pour les variétés compactes reliées par des anses fines
qu'on appliquera à la variété $M$ reliée par une anse à une sphère
(voir figure \ref{demo:anses}). 
\begin{figure}[h]
\begin{center}
\begin{picture}(0,0)%
\includegraphics{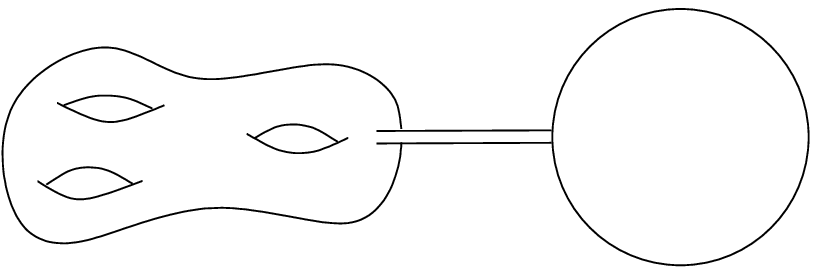}%
\end{picture}%
\setlength{\unitlength}{4144sp}%
\begingroup\makeatletter\ifx\SetFigFont\undefined%
\gdef\SetFigFont#1#2#3#4#5{%
  \reset@font\fontsize{#1}{#2pt}%
  \fontfamily{#3}\fontseries{#4}\fontshape{#5}%
  \selectfont}%
\fi\endgroup%
\begin{picture}(3706,1186)(1518,-1532)
\end{picture}%
\end{center}
\caption{\label{demo:anses}}
\end{figure}

Plus précisément, on considère une métrique $g\in C$
sur $M$ telle que $\mu_{\left[\frac n2\right],1}(M,g)>
\mu_k^c(S^n,C_\textrm{can})$, et pour tout réel $\delta>0$ et tout entier
$k\geq1$, on se donne une métrique $g_{k,\delta}$ sur la sphère $S^n$
conforme à la métrique canonique telle que 
$\mu_{\left[\frac n2\right],k}(S^n,g_{k,\delta})\leq
\mu_k^c(S^n,C_\textrm{can})+\delta$ et $\Vol(S^n,g_{k,\delta})=1$.
On fixe deux points $x\in M$ et $y\in S^n$ qui seront les points d'attache 
de l'anse, et on déforme les métriques
$g$ et $g_{k,\delta}$ en des métriques $g_\delta$ et $g_{k,\delta}'$
telles que $g_\delta$ (resp. $g_{k,\delta}'$) soit euclidienne sur
une boule de centre $x$ (resp. $y$) et de rayon $\delta$ et que 
$g_{k,\delta}'$ soit conforme à $g_{k,\delta}$ (c'est possible puisque 
la métrique canonique de la sphère est conformément plate).
En suivant un procédé déjà utilisé dans \cite{ces06} et \cite{ja06}, 
on déforme ensuite la métrique sur la boule euclidienne de centre $x$ afin 
d'obtenir la figure \ref{demo:anses}: pour tout $0<\varepsilon\leq\delta$, on 
définit une fonction $h_{1,\varepsilon}$ sur $M$ par $h_{1,\varepsilon}=1$ 
en dehors de la boule de centre $x$, et 
\begin{equation}\label{demo:cyl}
h_{1,\varepsilon}(r)=\left\{\begin{array}{ll}
e^{\frac L\varepsilon}&\textrm{si }0\leq r\leq 
\varepsilon e^{-\frac L\varepsilon},\\
\frac\varepsilon r&\textrm{si }\varepsilon e^{-\frac L\varepsilon}
\leq r\leq\varepsilon,\\
1 &\textrm{si }\varepsilon\leq r\leq\delta,
\end{array}\right.
\end{equation}
où $r$ est la coordonnée radiale sur la boule $B(x,\delta)$. Si on munit
$M$ de la métrique $h_{1,\varepsilon}^2g_\delta$, la boule $B(x,\varepsilon)$ 
devient isométrique à la réunion d'un cylindre de
longueur~$L$ et d'une boule euclidienne de rayon~$\varepsilon$ reposant sur 
le bord du cylindre. On peut projeter stéréographiquement ---~donc de manière 
conforme~--- cette boule sur une calotte sphérique reposant sur le
bord du cylindre et donc trouver une fonction
$h_{2,\varepsilon}$ telle que $(M,h_{2,\varepsilon}^2g_\delta)$ soit la 
variété de la figure \ref{demo:anses}. Comme la métrique $g_{k,\delta}'$ 
sur la sphère est conforme à la métrique canonique, on peut trouver 
une fonction $h_{3,\varepsilon}$ telle que
$(M,h_{3,\varepsilon}^2g_\delta)$ soit la variété obtenue en reliant 
$(M,g_\delta)$ et $(S^n,g_{k,\delta}')$ par une anse de longueur $L$ et de 
rayon $\varepsilon$ reposant sur les bords des boules de centre $x$ et $y$
et de rayon $\varepsilon$. Le théorème B de \cite{ac95} nous dis alors que
$(\mu_{\left[\frac n2\right],i}(M,h_{3,\varepsilon}^2g_\delta))_{i\geq1}$ 
tend vers la réunion de $(\mu_{\left[\frac n2\right],i}(M,g_\delta))_{i\geq1}$
et $(\mu_{\left[\frac n2\right],i}(S^n,g_{k,\delta}'))_{i\geq1}$ 
quand $\varepsilon\to0$.

On peut choisir les métriques $g_\delta$ et $g_{k,\delta}'$ telles que
$\tau(\delta)^{-1}g_\delta\leq g\leq\tau(\delta)g_\delta$ et
$\tau(\delta)^{-1}g_{k,\delta}'\leq g_{k,\delta}\leq\tau(\delta)
g_{k,\delta}'$ avec $\tau(\delta)\to0$ quand $\delta\to0$. Un résultat
de J.~Dodziuk (\cite{do82}, théorème 3.3) permet alors d'affirmer que 
$\mu_{\left[\frac n2\right],i}(M,g_\delta)\to\mu_{\left[\frac n2\right],i}
(M,g)$ et $\mu_{\left[\frac n2\right],i}(S^n,g_{k,\delta}')\to
\mu_i^c(S^n,C_\textrm{can})$ pour tout $i\geq1$ quand $\delta\to0$. 
On peut donc trouver une suite de fonctions $(h_\delta)$ telle que
$\mu_{\left[\frac n2\right],k}(M,h_\delta^2g_\delta)\to
\mu_k^c(S^n,C_\textrm{can})$ quand $\delta\to0$ (on utilise le fait
que le spectre de $(M,g)$ est minoré par $\mu_k^c(S^n,C_\textrm{can})$). 
Les métriques 
$h_\delta^2g_\delta$ ne sont pas conformes à la métrique $g$, 
mais on a $\tau(\delta)^{-1}h_\delta^2g_\delta\leq 
h_\delta^2g\leq\tau(\delta)h_\delta^2g_\delta$
ce qui assure que $\mu_{\left[\frac n2\right],k}(M,h_\delta^2g)$ tend
aussi vers $\mu_k^c(S^n,C_\textrm{can})$ quand $\delta$ tend vers zéro.
Comme la métrique $g$ peut être choisie de volume arbitrairement petit,
on peut en outre faire tendre le volume vers~1. On a donc bien 
$\mu_k^c(M,C)\leq\mu_k^c(S^n,C_\textrm{can})$.

On peut démontrer l'inégalité 
$\mu_k^c(S^n,C_\textrm{can})\leq k^{\frac2n}\mu_1^c(S^n,C_\textrm{can})$
en utilisant la même technique.
On se donne une métrique $g_\delta\in C_\textrm{can}$ sur $S^n$ telle que
$\mu_{\left[\frac n2\right],1}(S^n,g_\delta)\leq\mu_1^c(S^n,C_\textrm{can})
+\delta$
et $\Vol(S^n,g_\delta)=1$, et on déforme conformément la sphère 
---~en répétant $k-1$ fois le procédé précédent~--- de manière à
la rendre isométrique à $k$ sphères munies de la métrique $g$ attachées par 
des anses de rayon $\varepsilon$ (figure \ref{demo:sphere}).
\begin{figure}[h]
\begin{center}
\begin{picture}(0,0)%
\includegraphics{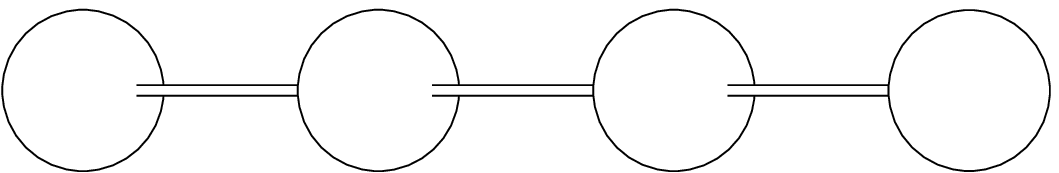}%
\end{picture}%
\setlength{\unitlength}{4144sp}%
\begingroup\makeatletter\ifx\SetFigFont\undefined%
\gdef\SetFigFont#1#2#3#4#5{%
  \reset@font\fontsize{#1}{#2pt}%
  \fontfamily{#3}\fontseries{#4}\fontshape{#5}%
  \selectfont}%
\fi\endgroup%
\begin{picture}(4808,756)(1546,-1236)
\end{picture}%
\end{center}
\caption{\label{demo:sphere}}
\end{figure}
Quand $\varepsilon$ tend
vers zéro, le volume de la sphère tend vers $k$ et les $k$ premières
valeurs propres sur les $\left[\frac n2\right]$-formes coexactes
tendent vers $\mu_{\left[\frac n2\right],1}(S^n,g_\delta)$, ce qui
permet d'affirmer que $\mu_k^c(S^n,C_\textrm{can})\leq
k^{\frac2n}\mu_{\left[\frac n2\right],1}(S^n,g_\delta)$. On conclut en faisant
tendre $\delta$ vers zéro.
\end{demo}

\noindent Pierre \textsc{Jammes}\\
Université d'Avignon\\
laboratoire de mathématiques\\
33 rue Louis Pasteur\\
F-84000 Avignon\\
\texttt{Pierre.Jammes@univ-avignon.fr}
\end{document}